\documentclass[11pt]{article}
\usepackage{float}
\usepackage{pifont}
\usepackage{rotating}
\usepackage[normalem]{ulem}
\usepackage[numbers]{natbib}
\usepackage[table]{xcolor}
\usepackage{tikz}
\usepackage{amsmath,amssymb,amsmath,amsthm,epsfig,graphicx,natbib}
\usepackage{color}
\usepackage{appendix}
\usepackage{enumitem}
\usepackage{authblk}
\usepackage{lscape}
\usepackage{subcaption}
\usepackage{bm}
\usepackage{lineno}
\usepackage{algorithm}
\usepackage{algorithmicx}
\usepackage{algpseudocode}
\usepackage{mwe}
\usepackage{colortbl}
\usepackage{multirow}
\usepackage{appendix}
\usepackage{soul}
\usepackage{caption}
\usepackage{pdfpages}

\usepackage{hyperref}
\hypersetup{%
  colorlinks=true,
  linkcolor=blue,
  citecolor=blue,
}
\newcommand{\RE}{{\mathbb R}}

\newtheorem{prop}{Proposition}

\newtheorem{lem}{Lemma}
\newtheorem{thm}{Theorem}
\newtheorem{ass}{Assumption}
\newtheorem{axi}{Axiom}
\newtheorem{ex}{Example}
\newtheorem{counterex}{Counterexample}
\newtheorem{rem}{Remark}
\newtheorem*{ack}{Acknowledgement}

\newcommand{\Keywords}[1]{\par\noindent
{\small{\bf Keywords\/}: #1}}

\begin{document}

\includepdf{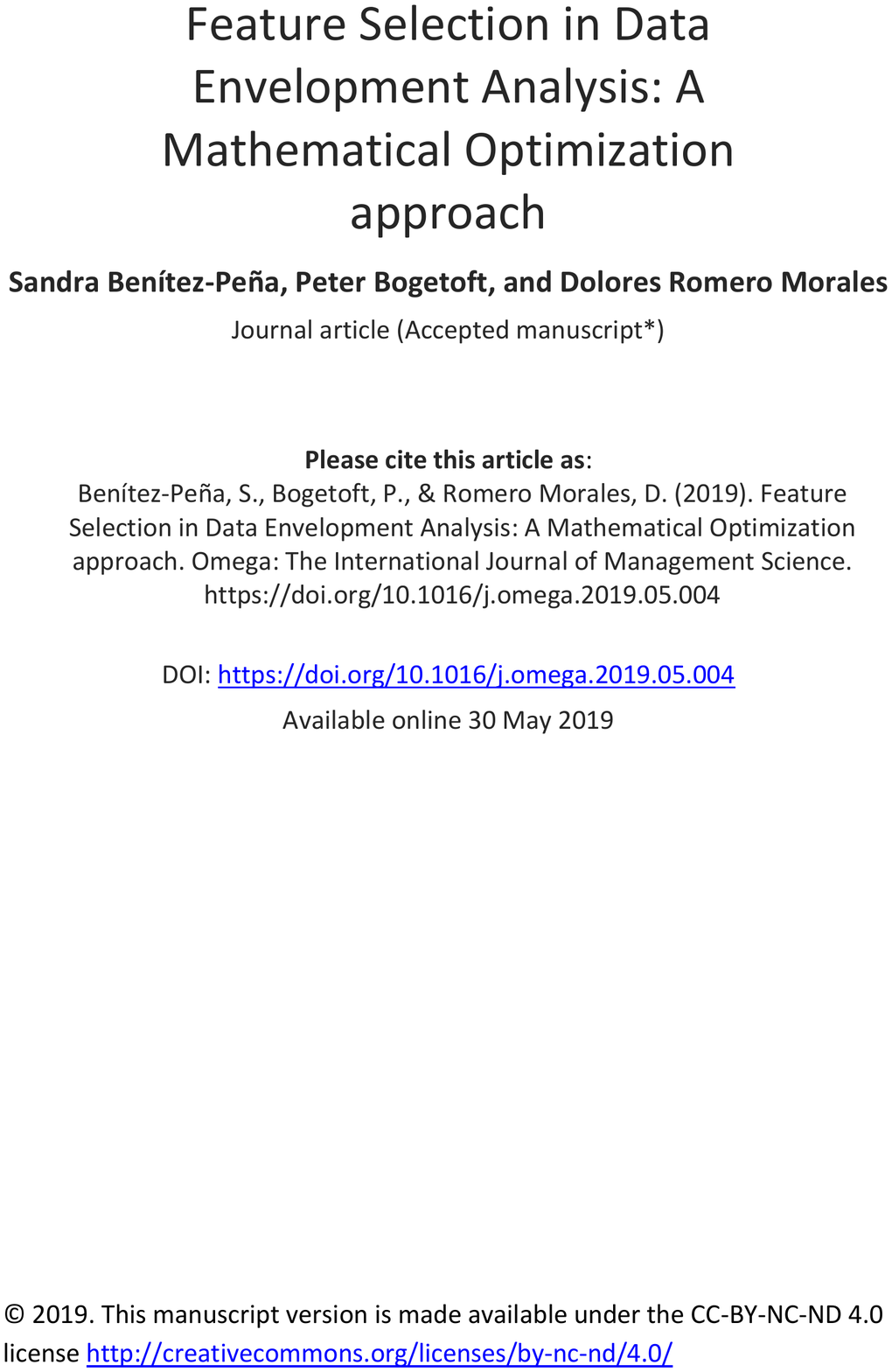}

\title{Feature Selection in Data Envelopment Analysis: A Mathematical Optimization approach}

\author[1]{Sandra Ben\'{\i}tez Pe\~na}
\author[2]{Peter Bogetoft}
\author[2]{Dolores Romero Morales}
\affil[1]{Instituto de Matem\'aticas de la Universidad de Sevilla (IMUS), Seville, Spain \newline
        \tt{sbenitez1@us.es}}
\affil[2]{Copenhagen Business School, Frederiksberg, Denmark \newline
        \tt{\{pb.eco,drm.eco\}@cbs.dk}}

\date{\today}

\maketitle

\begin{abstract}
\noindent
This paper proposes an integrative approach to feature (input and output) selection in Data Envelopment Analysis (DEA). The DEA model is enriched with zero-one decision variables modelling the selection of features, yielding a Mixed Integer Linear Programming formulation. This single-model approach can handle different objective functions as well as constraints to incorporate desirable properties from the real-world application. Our approach is illustrated on the benchmarking of electricity Distribution System Operators (DSOs). The numerical results highlight the advantages of our single-model approach provide to the user, in terms of making the choice of the number of features, as well as modeling their costs and their nature.
\vspace*{1cm}

\noindent
\end{abstract}

\newpage

\section{Introduction \label{sec:intro}}

Organisations need to know whether they are using the best practices to produce their products and services, and to do so they benchmark their performance with that of others. There are many documented examples of the use of benchmarking in the literature from both the private and the public sector, such as, airlines, banks, hospitals, universities, manufacturers, schools, and municipalities, see \cite{bogetoft2013performance}, and references therein.

Within benchmarking, Data Envelopment Analysis (DEA) \citep{CHARNES1978429} is one of the most widely used tools, \cite{COOK201945,emrouznejad2018survey,jiang2015dearank,landete2017robust,li2017dynamic,petersen2018directional,RUIZ20161,RUIZ2018}. It aims at benchmarking the performance of decision marking units (DMUs), which use the same types of inputs and produce the same types of outputs, against each other. DEA calculates an efficiency score for each of the DMUs, so that DMUs with a score equal to one are in the so-called efficient frontier. DMUs outside the efficient frontier are deemed as underperforming, and a further analysis gives insights as to what they can do to improve their efficiency. The efficiency of DMUs in DEA is measured as the weighted summation of the outputs divided by the weighted summation of the inputs, and the weights are found solving a Linear Programming problem for each DMU. DEA model specification, in the form of feature (where the term feature is used to refer to either outputs, inputs or environmental variables) selection, has a significant impact on the shape of the efficient frontier in DEA as well as the insights given to the inefficient DMUs \cite{golany1989application}. Moreover, it is known to improve the discriminatory power of DEA models \cite{bogetoft2010benchmarking}. Our paper proposes and investigates a mathematical optimization approach for feature selection in DEA.

In benchmarking projects, as in most applied statistical analyses, one of the most challenging tasks is the choice of the DEA model specification. First, a good model should make conceptual sense not only from the theoretical but also from a practical point of view. The interpretation must be easy to understand and the properties of the model must be natural. This contributes to the acceptance of the model by stakeholders and provides a safeguard against spurious models developed without much understanding of the industry. More precisely, this has to do with the choice of outputs in DEA that are natural cost drivers and with functional forms that, for example, have reasonable returns to scale and curvature properties. Second, it is important to guide the search for a good model with classical statistical tests. We typically seek models that have significant features with the right signs and that do not leave a large unexplained variation. Third, intuition and experience is a less stringent but important safeguard against false model specifications and the over- or underuse of data to draw false conclusions. It is important that the models produce results that are not that different from the results one would have found in other data sets, e.g., from other countries or related industries. The intuition and experience must  be used with caution. We may screen away extraordinary but true results (Type 1 error) and we may go for a more common set of results based on false models (Type 2 error). One aspect of this is that one will tend to be more confident in a specification of inputs and outputs that leads to comparable results in alternative estimation approaches, e.g., in the DEA and Stochastic Frontier Analysis models. Finally, the choice of model specification has to be pragmatic. We need to take into account the availability of data as well as what the model is going to be used for. In benchmarking, it matters if the model is used to learn best practices, to reallocate resources between entities or to directly incentivize firms or managers by performance based payment schemes. Our approach gives a tool that can support the selection of features in benchmarking, allowing the user to navigate through a large amount of DEA models and a large amount of constraints modeling knowledge in the form of intuition and experience, in an efficient manner.

The complexity of the model specification phase partially explains the lack of enough guidance in the literature at this respect, \cite{cook2014data,luo2012input,SOLEIMANIDAMANEH20095146}, and most of the effort goes into the analysis and interpretation of a given DEA model. With the strand of literature on feature selection, the most common approach is to use a priori rules based on Statistical Analysis (such as correlations, dimensionality reduction techniques, and regression), and Information Theory (such as AIC or Shannon entropy). Alternatively, an ex-post analysis of the sensitivity of the efficient frontier to additional features can be run to detect whether relevant features have been left out. See \cite{adler2010improving,fernandez2018stepwise,li2017variable,nataraja2011guidelines,pastor2002statistical,sirvent2005monte,SOLEIMANIDAMANEH20095146,wagner2007stepwise}, and references therein. Recently, there have been attempts to use LASSO techniques from Statistical Learning to build sparse benchmarking models, i.e., models using just a few features, \cite{caiLASSODEA,LEE2018,qin2014joint}.

In this paper, the DEA Linear Programming formulation is enriched with zero-one decision variables modelling the selection of features for different objective functions, such as the average efficiency or the squared distance to the ideal point where all DMUs are efficient, and for different set of constraints that incorporate knowledge from the industrial application, such as bounds on the weights as well as costs on the features. This yields either a Mixed Integer Linear Programming (MILP) problem, or a Mixed Integer Quadratic Programming (MIQP) one. Thus, in contrast to the existing literature, that tends to combine statistical analysis with the mathematical programming based DEA, we propose an approach that is entirely driven by mathematical optimization.
We illustrate our models in the benchmarking of electricity Distribution System Operators (DSOs), where there is a pool of 100 potential outputs.

The contributions of our approach are threefold. First, our single-model mathematical approach can guide better the selection of features: it controls directly the number of chosen features, as opposed to techniques based on seeking sparsity, being thus able to quantify the added value of additional features; works directly with the original features, as opposed to dimensionality reduction techniques, which create artificial features that are difficult to interpret; and can derive a collection of models by shaping in alternative ways the distribution of the efficiencies, using different objective functions that focus on different groups of DMUs, which can be combined through, for instance, Shannon entropy \cite{SOLEIMANIDAMANEH20095146}. Second, while the previous literature has focused on the choice of variables from a small set of candidates, e.g., \cite{LEE2018}, in the era of Big Data, the set of alternatives to choose from is expanding at a fast pace, and the challenge is often not the lack of data, but the abundance of data, \cite{ZHU2018291}. In the numerical section, we show how our MILP/MIQP approach is able to make the selection from a large pool of outputs. Third, we introduce an element of game theory when selecting features.
In applied projects, the evaluated DMUs will typically try to influence the feature selection since this will affect how one firm is evaluated relative to others. It is therefore important to think about the conflict the DMUs (the players of the game) when choosing the set of features (the strategies of the game) used in the calculation of the efficiencies (outcome). We illustrate the results for a simpler game setting where the strategies are derived from the individual and the joint models.

The reminder of the paper is structured as follows. In Section \ref{sec:OS}, we introduce the individual feature selection problem where the selection is tailored to a given DMU. In Section \ref{sec:singleOS}, we introduce the joint feature selection problem where the selection is imposed to be the same for all DMUs. Section \ref{sec:numerical} is devoted to the illustration of our models in the benchmarking of electricity Distribution System Operators (DSOs). We end the paper in Section \ref{sec:conclusions} with some conclusions and lines for future research.

\section{The individual selection model\label{sec:OS}}

In this section, and for an individual DMU, we propose a Mixed Integer Linear Programming (MILP) formulation to select outputs in Data Envelopment Analysis (DEA).

Consider $K$ DMUs (indexed by $k$), using $I$ inputs (indexed by $i$), and producing $O$ outputs (indexed by $o$). DMU $k$ uses vector of inputs $\mathbf{x}^{(k)}\in \RE_+^I$ to produce vector of outputs $\mathbf{y}^{(k)}\in \RE_+^O$. Let $E^{(k)}$ be the so-called Farrell input-oriented efficiency of DMU $k$, which is the optimal solution value to a DEA model. Our goal is to select the $p$ outputs from the $O$ potential ones that yield the maximal efficiency for DMU $k$. We first start with the classical formulation of the problem which solves a Linear Programming model, and subsequently include the output selection decision variables. The output selection model for DMU $k$ is enriched with input selection decision variables, as well as constraints modeling desirable properties about the selected features. Please note that our approach can easily be extended to the use of other efficiency measures, including the output-based Farrell efficiency.

\subsection{The selection model for a DMU}

We start with the formulation of the classical DEA model, in which we can make use of the $O$ outputs available. The input-oriented efficiency of DMU $k$, $E^{(k)}$, in a DEA model with constant returns to scale (CRS) is equal to the optimal solution value of the following Linear Programming formulation,
\begin{eqnarray}
E^{(k)} = \mbox{maximize}_{(\bm{\alpha}^{(k)},\bm{\beta}^{(k)})} \ \sum_{o=1}^O \beta_o^{(k)} y_o^{(k)} \label{eq:of}
\end{eqnarray}
\noindent {\em s.t.} \hfill (DEA$^{(k)}$)
 \begin{eqnarray}
  \displaystyle
   \sum_{o=1}^O \beta_o^{(k)} y_o^{(j)} - \sum_{i=1}^I \alpha_i^{(k)} x_i^{(j)} \le 0 && \forall  j=1,\ldots,K\label{eq:DEA1}\\
   \sum_{i=1}^I \alpha_i^{(k)} x_i^{(k)} =1 && \label{eq:DEA2}\\
   \alpha^{(k)} \in \RE^I_+ && \label{eq:DEA3}\\
   \beta^{(k)} \in \RE^O_+, && \label{eq:DEA4}
 \end{eqnarray}
where $\alpha_i^{(k)}$ is the weight for input $i$ and $\beta_o^{(k)}$ the weight for output $o$. (DEA$^{(k)}$) has $K+1$ linear constraints and $I+ O$ continuous variables, and thus can be solved efficiently even for large problem instances.

We continue with the model in which $p$ outputs are to be selected from the $O$ available ones such that the efficiency of DMU $k$ is maximized. Let $z^{(k)}_o$ be equal to 1 if output $o$ can be used in the calculation of the efficiency  of DMU $k$, and 0 otherwise. Let $E^{(k)}(\mathbf{z}^{(k)})$ denote the corresponding efficiency. The decision variables $\beta_o^{(k)}$ and $\alpha_i^{(k)}$ are defined as above. The Output Selection for DMU $k$, where $p$ outputs must be selected such that $E^{(k)}(\mathbf{z}^{(k)})$ is maximized, can be written as the following MILP:
\begin{eqnarray}
v^{(k)}(p) := \mbox{maximize}_{(\bm{\alpha}^{(k)},\bm{\beta}^{(k)},\mathbf{z}^{(k)})} \ \displaystyle \sum_{o=1}^O \beta_o^{(k)} y_o^{(k)}
\end{eqnarray}
\noindent {\em s.t.} \hfill (OSDEA$^{(k)}(p)$)
 \begin{eqnarray}
  \displaystyle
  \eqref{eq:DEA1}-\eqref{eq:DEA4} && \nonumber\\
   \beta^{(k)}_o \le M z^{(k)}_o && \forall o=1,\ldots,O\label{eq:DEA5OS}\\
   \sum_{o=1}^O z^{(k)}_o = p && \label{eq:DEA6OS}\\
   z^{(k)}_o \in \{0,1\} && \forall o=1,\ldots,O\label{eq:DEA7OS},
 \end{eqnarray}
where $M$ is a big constant. Constraints \eqref{eq:DEA1}-\eqref{eq:DEA4} were already present in the classical DEA model. Constraints (\ref{eq:DEA5OS}) make sure that the selection variables $z^{(k)}_o$ are well defined: if $z^{(k)}_o$ equals $0$, then $\beta_o^{(k)}$ equals $0$ too. Constraint (\ref{eq:DEA6OS}) models the number of features to be selected. Finally, constraints (\ref{eq:DEA7OS}) relate to the range of decision variables $z^{(k)}_o$. (OSDEA$^{(k)}(p)$) has $K+O+2$ linear constraints and $I+ 2 \, O$ variables, where $I+O$ are continuous and $O$ are binary ones. Our numerical experiments show that this problem can be solved efficiently, although the solution time is affected by the value of the $M$ constant. The value of $M$, and thus the computational burden of the problem, can be reduced using an upper bound on the weight associated with output $o$, for each $o=1,\ldots,O$. It is not difficult to see that, without loss of optimality,
$\beta^{(k)}_o=0$ if $y_o^{(k)}=0$, and thus $z^{(k)}_o=0$. Otherwise, $\beta^{(k)}_o \le \frac{1}{y_o^{(k)}}$, by combining \eqref{eq:DEA1} and \eqref{eq:DEA2}. Thus, constraints \eqref{eq:DEA5OS} can be tighten to
 \begin{eqnarray*}
  \displaystyle
   \beta^{(k)}_o = 0 && \forall o=1,\ldots,O, \mbox{such that } y_o^{(k)}=0\\
   \beta^{(k)}_o \le \frac{1}{y_o^{(k)}} \,\, z^{(k)}_o  && \forall o=1,\ldots,O, \mbox{such that } y_o^{(k)}>0.
 \end{eqnarray*}

Let $z^{(k)}(p)$ denote the optimal selection variables to (OSDEA$^{(k)}(p)$), i.e., the $p$ outputs that yield the maximum efficiency for DMU $k$. Thus, the optimal solution value to (OSDEA$^{(k)}(p)$), denoted above by $v^{(k)}(p)$, is equal to $E^{(k)}(z^{(k)}(p))$. A few things are known about the maximum efficiency $v^{(k)}(p)$ as function of $p$. The efficiency $v^{(k)}(p)$ is non decreasing in $p$, i.e., the more outputs we select the better the efficiency of DMU $k$ can be. Moreover, in the extreme case when all outputs are considered, we have that $v^{(k)}(O)=E^{(k)}$. Thus, a plausible strategy to choose the value of $p$ is to look at the marginal contribution of an additional feature, i.e., $v^{(k)}(p+1) - v^{(k)}(p)$, and stop when this is below a threshold.

\subsection{Extensions}
\label{sec:extensions}

In this section we discuss several interesting extensions that can be carried out using the previous model as a basis. First, we present the formulation when both inputs and outputs are to be selected, all at once. Second, we model constraints on the weights attached to the outputs. Finally, we discuss how other attributes attached to the outputs, such as costs and correlations, may constrain the feature selection.

\subsubsection{Selection of inputs and outputs}

Note that up to now, and for the sake of clarity, we have focused on the selection of outputs. The selection of $\tilde{p}$ inputs from the $I$ potential ones can be included in a similar fashion. Indeed, let us consider the new binary variables $\tilde{z}^{(k)}_i$, equal to 1 if input $i$ can be used in the calculation of the efficiency for DMU $k$, and 0 otherwise. Hence, the Feature Selection for DMU $k$, (FSDEA$^{(k)}(p)$), where $\tilde{p}$ inputs and $p$ outputs are selected, can be written also as an MILP
\begin{eqnarray}
\mbox{maximize}_{(\bm{\alpha},\bm{\beta},\mathbf{z}^{(k)},\tilde{\mathbf{z}}^{(k)})} \ \displaystyle \sum_{o=1}^O \beta_o^{(k)} y_o^{(k)}
\end{eqnarray}
\noindent {\em s.t.} \hfill (FSDEA$^{(k)}(p)$)
 \begin{eqnarray}
  \displaystyle
  \eqref{eq:DEA1}-\eqref{eq:DEA7OS}&& \nonumber\\
   \alpha^{(k)}_i \le \tilde{M} \tilde{z}^{(k)}_i && \forall i=1,\ldots,I\label{eq:DEA252OS}\\
   \sum_{i=1}^I \tilde{z}^{(k)}_i = \tilde{p} && \label{eq:DEA262OS}\\
   \tilde{z}^{(k)}_i \in \{0,1\} && \forall i=1,\ldots,I\label{eq:DEA272OS},
 \end{eqnarray}
where $ \tilde{M}$ is another big constant. Constraints (\ref{eq:DEA252OS})--(\ref{eq:DEA272OS}) are the counterparts of (\ref{eq:DEA5OS})--(\ref{eq:DEA7OS}) but modelling the selection of inputs instead of outputs. (FSDEA$^{(k)}(p)$) has $K+O+I+3$ linear constraints and $2 \, I+ 2 \, O$ variables, where half of them are continuous and the other half binary. Although running times are not an issue for this model, one can lower them even further by finding tighter values of $M$ and $ \tilde{M}$. As above, this can be done using bounds on the inputs and the outputs.

\subsubsection{Modeling constraints on weights}

Our (OSDEA$^{(k)}(p)$) improves the discriminatory power of the DEA model by focusing on a few outputs, and eliminating the rest from the calculation of the efficiency of DMU $k$. There is a strand of literature that, using also as a basis the discriminatory power, argue the necessity of controlling the values of the weights \citep{allen1997weights,GREEN1996461,joro2015data,PODINOVSKI2016916,ramon2010choice,doi101002ev1441}. In these works, it is assumed that we have upper and lower bounds on the weight $\beta^{(k)}_o$, say, $L^{(k)}_{o}$ and $U^{(k)}_{o}$, for $o = 1,\ldots,O$, i.e.,
  \begin{eqnarray}
  \displaystyle
   L^{(k)}_{o} \leq \beta^{(k)}_o \leq U^{(k)}_{o} & \forall o = 1,\ldots,O. \label{eq:lowerupperbound}
 \end{eqnarray}
Gathering non trivial values for these bounds is not a straightforward task for the user in the presence of many outputs, as in dataset on benchmarking of electricity DSOs in our numerical section. In any case, we can enrich our (OSDEA$^{(k)}(p)$), to not only control whether an output can be used, but also the range of values for the corresponding weight. These bounds can be incorporated in constraints (\ref{eq:DEA5OS}) in (OSDEA$^{(k)}(p)$) yielding
$$
  \begin{array}{ll}
  \displaystyle
   L^{(k)}_{o}z^{(k)}_o \leq \beta^{(k)}_o \leq U^{(k)}_{o} z^{(k)}_o & \forall o = 1,\ldots,O.
 \end{array}
 \eqno{(7')}
$$
There are a few observations to be made. First, the knowledge of upper bounds on the weights naturally tightens the value of $M$. Second, if there are meaningful lower bounds on the weights, i.e., if $L^{(k)}_{o}> 0$, then $z^{(k)}_o$ must be equal to $1$. Third, these positive lower bounds make the
selection problem (OSDEA$^{(k)}(p)$) infeasible for small values of $p$. Indeed, this is the case when there are more than $p$ outputs with a positive lower bound.

\subsubsection{Modeling attributes of the outputs}

Outputs may have attributes attached to them, which may affect the selection. We will model two of those.

First, we will consider that outputs are different in nature and therefore we will attach a different cost to them. Let $c_o$ denote the cost associated with output $y_o, o = 1, \ldots, O$, which can measure the collection and the verification of this output in a repeated setting. To select $p$ outputs so that their total cost does not exceed a given amount $C$, we need to add to (OSDEA$^{(k)}(p)$) the following constraint
 \begin{eqnarray}
  \displaystyle
    \sum_{o=1}^O c_o z^{(k)}_o \leq C. && \label{eq:DEAOScost}
 \end{eqnarray}

Second, we can consider the outputs being partitioned into $S$ clusters, with outputs within a cluster being similar in terms of what they measure. In the context of benchmarking electricity Distribution System Operators (network companies), clusters may related to the many different measurements of connections, transformers, lines, cables, etc. Let $\mathcal{H}= \{H_1, \ldots, H_S\}$ denote the partitioning of the outputs, namely $H_\ell \cap H_s = \emptyset$ and $\cup_{\ell=1}^S H_\ell=\{1,\ldots,O\}$. Given the similarity of outputs within a cluster, we will impose that at most (respectively, at least) $p^{\rm (max)}_\ell$ (respectively, $p^{\rm (min)}_\ell$) outputs can be selected from $H_\ell$. In order to do so, we need to add to (OSDEA$^{(k)}(p)$) the following constraint
 \begin{eqnarray}
  \displaystyle
    \sum_{o \in H_\ell} z^{(k)}_o \leq p^{\rm (max)}_\ell && \forall \ell=1,\ldots,S. \label{eq:DEAOScluster1}\\
    \sum_{o \in H_\ell} z^{(k)}_o \geq p^{\rm (min)}_\ell && \forall \ell=1,\ldots,S. \label{eq:DEAOScluster2}
 \end{eqnarray}

Finally, we have correlation $\rho_{oo'}$ between outputs $o$ and $o'$ as another attribute. If two outputs are highly correlated, we may be interested in using only one of them, since the information they provide is almost the same and can derive in the problem of multicollinearity \cite{Bertsimas2016ORF}. Hence, we want to impose that if $\rho_{oo'}$ is greater than a preselected threshold, then outputs $o$ and $o'$ cannot be chosen simultaneously. We can model this by first defining a 0--1 matrix $R$, in which $R_{oo'}=0$ if $\rho_{oo'}$ is lower than the threshold, and $1$ otherwise. Then, we have simply to add to (OSDEA$^{(k)}(p)$) the constraints
 \begin{eqnarray}
  \displaystyle
z^{(k)}_o + z^{(k)}_{o'} \leq 2- R_{oo'}, & \forall o < o'.\label{eq:correlmodel}
\end{eqnarray}
The choice of the threshold have to be done with care, since some works like \cite{nunamaker1985using} suggest that the addition of a highly correlated variable may increase the efficiency.

Throughout this section, we have made the selection of outputs individually for DMU $k$ with the goal to maximize the efficiency of DMU $k$. Therefore, for two different DMUs, $k$ and $ k'$, the selected outputs, $z^{(k)}(p)$ and $z^{(k')}(p)$, may differ. In model specification one is interested in finding the most discriminatory features in order to build a valid model for all DMUs. With this in mind, we propose in the next section a mathematical optimization model that selects the outputs jointly for all DMUs, ensuring they will be the same ones for all DMUs.

\section{The joint selection model\label{sec:singleOS}}

In this section, we address the problem in which the selected outputs have to be the same for all DMUs. First, this joint selection is made maximizing the average efficiency of all DMUs, yielding an MILP formulation. The model can be enriched as in previous section with input selection decision variables, as well as constraints modeling desirable properties about the selected features. Second, we propose alternatives to the maximization of the average efficiency when making the joint selection of outputs, such as the maximization of the weighted average efficiency, the minimum efficiency, or a percentile of the efficiencies. The joint selection model can again be formulated as an MILP problem. We also consider the minimization of the square of the Euclidean distance from each DMU efficiency to the ideal value of 1, where the joint selection model can be rewritten as a Mixed Integer Quadratic Programming problem.

\subsection{The selection model for all DMUs}

To obtain all efficiencies $E^{(k)}$ simultaneously, one can solve the following single-objective Linear Programming formulation
\begin{eqnarray}
\frac{1}{K} \sum_{k=1}^K E^{(k)} = \mbox{maximize}_{(\bm{\alpha},\bm{\beta})} \ \displaystyle \frac{1}{K} \sum_{k=1}^K \sum_{o=1}^O   \beta_o^{(k)}
y_o^{(k)} \label{eqn:ofsingle}
\end{eqnarray}
\noindent {\em s.t.}
 \begin{eqnarray}
  \displaystyle
   \sum_{o=1}^O \beta_o^{(k)} y_o^{(j)} - \sum_{i=1}^I \alpha_i^{(k)} x_i^{(j)} \le 0 && \forall  j=1,\ldots,K; \forall k=1,\ldots,K \label{eq:DEA1single}\\
   \sum_{i=1}^I \alpha_i^{(k)} x_i^{(k)} =1 && \forall k=1,\ldots,K \label{eq:DEA2single}\\
   \alpha \in \RE^{I \cdot K}_+ && \label{eq:DEA3single}\\
   \beta \in \RE^{O \cdot K}_+. && \label{eq:DEA4single}
 \end{eqnarray}
It is easy to see that this problem decomposes by DMU, and that each of the subproblems are equivalent to (DEA$^{(k)}$), which optimal solution value is $E^{(k)}$.

We continue with the model in which $p$ outputs are to be selected from the $O$ available ones, the same ones for all DMUs. The goal in this section is to maximize the average efficiency across all DMUs. Let $z_o$ be equal to 1 if output $o$ can be used in the calculation of the efficiencies, and 0 otherwise. The decision variables $\beta_o^{(k)}$ and $\alpha_i^{(k)}$ are defined as above. The Output Selection for DEA problem, (OSDEA$(p)$), where $p$ outputs must be selected such that the average efficiency across all DMUs is maximized, can be written as the following MILP:
\begin{eqnarray}
v(p) := \mbox{maximize}_{(\bm{\alpha},\bm{\beta},\mathbf{z})} \ \displaystyle \frac{1}{K} \sum_{k=1}^K \sum_{o=1}^O   \beta_o^{(k)} y_o^{(k)}
\end{eqnarray}
\noindent {\em s.t.} \hfill (OSDEA$(p)$)
 \begin{eqnarray}
  \displaystyle
  \eqref{eq:DEA1single}-\eqref{eq:DEA4single} && \nonumber\\
   \beta^{(k)}_o \le M z_o && \forall o=1,\ldots,O; \forall k=1,\ldots,K\label{eq:DEA5OSNEW}\\
   \sum_{o=1}^O z_o = p && \label{eq:DEA6OSNEW}\\
   z_o \in \{0,1\} && \forall o=1,\ldots,O\label{eq:DEA7OSNEW},
 \end{eqnarray}
where $M$ is a big constant. Constraints \eqref{eq:DEA1single}-\eqref{eq:DEA4single} are necessary to find the weights of the inputs and the outputs that the efficiency for each DMU. Constraints (\ref{eq:DEA5OSNEW}) make sure that the selection variables $z_o$ are well defined with respect to $\beta_o^{(k)}$. Constraint (\ref{eq:DEA6OSNEW}) models the number of features to be selected. Finally, constraints (\ref{eq:DEA7OSNEW}) relate to the range of decision variables $z_o$. (OSDEA$(p)$) has $K(K+O+1)+1$ linear constraints and $K(O+I)+O$ variables, where $K(O+I)$ are continuous and $O$ are binary ones. We have multiplied the size of the problem by $K$, except for the number of binary variables, which are still one per output. Our numerical experiments show that this problem can still be solved efficiently. Moreover, and as in previous section, the computational burden of the problem depends on the value $M$ and we can tighten it using similar bounds. As before, we might extend the model as in Section \ref{sec:extensions}, with input selection decision variables, as well as constraints to model desirable properties of the outputs.

Let $z(p)$ denote the optimal selection variables to (OSDEA$(p)$), i.e., the $p$ outputs that yield the maximum average efficiency. Thus, the optimal solution value to (OSDEA$(p)$), denoted above by $v(p)$, is equal to $\frac{1}{K} \sum_{k=1}^K E^{(k)}(z(p))$. In general, we have that $E^{(k)}(z(p)) \le E^{(k)}(z^{(k)}(p))$, since $z^{(k)}(p)$ is the best strategy for DMU $k$. As in previous section, the maximum average efficiency $v(p)$ is non decreasing in $p$, i.e., the more outputs we select the better the average efficiency can be. In the limit case, we have that $v(O)=\frac{1}{K} \sum_{k=1}^K E^{(k)}$. The number of selected outputs $p$ is a parameter of our model. The user should make the choice of $p$ after inspecting the curve $v(p)$. As before, a plausible strategy to choose the value of $p$ is to look at the marginal contribution of an additional feature, i.e., $v(p+1) - v(p)$, and stop when this is below a threshold. The question is whether this marginal contribution is nonincreasing, i.e., $v(p+2) - v(p+1) \le v(p+1) - v(p)$, for all $p=1,\ldots,O-2$. Below, we show a toy example where this inequality is not satisfied, and thus $v(p)$ is not a concave function of $p$. In the numerical section, devoted to the benchmarking of electricity DSOs, the function $v(\cdot)$ that we obtain empirically is concave, and thus, not convex.

\begin{counterex} \label{counterex:concave}
Consider $5$ DMUs, each one described by a single input and four different outputs, as can be seen in Table~\ref{tab:Counter2}. When performing the feature selection procedure, the results that we obtain are the following. In the case of selecting just one output, say $p=1$, the procedure chooses ``Output 1'' and the obtained efficiencies are then just the same as the values of ``Output 1" for each DMU. Hence, the average efficiency is 0.8. When two outputs are selected, the procedure chooses ``Output 1'' and ``Output 2''. These outputs make the average efficiency to be 0.867. Furthermore, if three outputs are selected, the procedure chooses ``Output 2'', ``Output 3'' and ``Output 4''. These outputs make all the DMUs efficient (i.e., efficiency equal to $1$) and thus the average efficiency is 1. Clearly,
\[v(3)-v(2) >v(2)-v(1),\]
and therefore $v(\cdot)$ is not concave.
\begin{table}
$$
\begin{tabular}{|c|c|cccc|}
   \hline
   DMU & Input 1 & Output 1 & Output 2 & Output 3 & Output 4 \\
   \hline
   1 & 1 & 0.6 & $\frac{1}{3}$ & $\frac{1}{3}$ & $\frac{1}{3}$ \\
   2 & 1 & 0.7 & $\frac{1}{3}$ & $\frac{1}{3}$ & $\frac{1}{3}$ \\
   3 & 1 & 0.8 & 1 & 0 & 0 \\
   4 & 1 & 0.9 & 0 & 1 & 0 \\
   5 & 1 & 1   & 0 & 0 & 1 \\
   \hline
\end{tabular}
$$
\caption{Toy example for which $v(\cdot)$ is not concave}
\end{table}\label{tab:Counter2}
\end{counterex}

A greedy approach is provided in \cite{pastor2002statistical} to address the feature selection problem in a nested fashion. In short, this greedy nested procedure works as follows.
For $p=1$, (OSDEA$(p)$) is solved to optimality. Let $o^{\rm}(1)$ be its best output. For $p=2$, (OSDEA$(p)$) is solved to optimality, with the additional constraint that $z_{o^{\rm}(1)}=1$. Let  $o^{\rm}(2)$ be its best output. In general, for $p$, (OSDEA$(p)$) is solved to optimality, with the additional constraints that $z_{o^{\rm}(1)}=z_{o^{\rm}(2)}=\ldots=z_{o^{\rm}(p-1)}=1$. Let $o^{\rm}(p)$ be its best output. Clearly, this greedy approach returns a sequence of outputs that is nested, i.e., the outputs selected in iteration $p-1$ will also be selected in iteration $p$, for all $p$.
The following is a toy example that illustrates that the approach in \cite{pastor2002statistical} does not provide, in general, the optimal solution to (OSDEA$(p)$).
\begin{counterex} \label{counterex:nested}
Consider $4$ DMUs, each one described by a single input and three different outputs, as can be seen in Table~\ref{tab:Counter}. When performing the feature selection procedure, the results that we obtain are the following. In the case of selecting just one output, say $p=1$, the procedure chooses ``Output 1'' and the obtained efficiencies are then just the same as the values of ``Output 1" for each DMU. When two outputs are selected, the procedure chooses ``Output 2'' and ``Output 3''. These outputs make all the DMUs efficient. However, if either ``Output 1'' and ``Output 2'' or ``Output 1'' and ``Output 3'' were used instead, the efficiencies would be \{0.85,0.9,0.95,1\} and \{1,1,0.927,1\}, respectively.
\begin{table}
$$
\begin{tabular}{|c|c|ccc|}
 \hline
   DMU & Input 1 & Output 1 & Output 2 & Output 3 \\
   \hline
   1 & 1 & 0.85 & 0.2 & 0.8 \\
   2 & 1 & 0.95 & 0.4 & 0.6 \\
   3 & 1 & 0.9 & 0.6 & 0.4 \\
   4 & 1 & 1 & 0.8 & 0.2 \\
\hline
\end{tabular}
$$
\caption{Toy example for which the approach in \cite{pastor2002statistical} does not provide the optimal solution to (OSDEA$(p)$)}
\end{table}\label{tab:Counter}
\end{counterex}

\subsection{Alternative objective functions}

In (OSDEA$(p)$), we maximize the average efficiency across all DMUs. In this section, we propose other objective functions $\phi()$ to select the outputs.

A straightforward generalization would be to consider the weighted average efficiency.
\begin{equation}
\label{eq:quadratic}
\phi^{\rm (w)}(\bm{\alpha},\bm{\beta},\mathbf{z}) = \frac{1}{K} \sum_{k=1}^K \sum_{o=1}^O \omega^{(k)} \beta_o^{(k)} y_o^{(k)}.
\end{equation}
This is relevant if the DMUs are not equally important. If there is only one input, say cost in $\$$, we could for example use $w^{(k)}=x^{(k)}$, and (\ref{eq:quadratic}) would correspond to minimizing the total sector loss from inefficiency.

Instead of the weighted average efficiency, one could be interested in measuring how far each DMU is from efficiency. This can be measured with the following quadratic function
\begin{equation}
\label{eq:quadratic2}
\phi^{\rm (q)}(\bm{\alpha},\bm{\beta},\mathbf{z}) = \frac{1}{K} \sum_{k=1}^K (1-  \sum_{o=1}^O \beta_o^{(k)}y_o^{(k)})^2.
\end{equation}
Alternatively, our goal could have been maximizing the worst efficiency, i.e., the minimum one
\begin{equation}
\label{eq:quadratic3}
\phi^{\rm (m)}(\bm{\alpha},\bm{\beta},\mathbf{z})= \min_{k=1,\ldots,K} \sum_{o=1}^O \beta_o^{(k)} y_o^{(k)}.
\end{equation}
This is relevant, for example, if outputs are selected with the aim of being Rawlsian fair towards all DMUs.
Instead of the minimum, we could have optimized another $\pi$-percentile, $\pi =1,\ldots,100$, of the efficiency distribution. Assuming that the efficiencies are given in non-decreasing order,
$ \sum_{o=1}^O\beta_o^{(k)}\mathbf{y}_o^{(k)} \le \sum_{o=1}^O\beta_o^{(k+1)}y_o^{(k+1)}$, for all $k$, we would have
\begin{equation}
\label{eq:percentile}
\phi^{(\pi)}(\bm{\alpha},\bm{\beta},\mathbf{z}) = \sum_{o=1}^O\beta_o^{(k(\pi))}y_o^{(k(\pi))},
\end{equation}
with $k(\pi)=\lfloor K \, \dfrac{\pi}{100}\rfloor.$

The Output Selection for DEA problem where the goal is to maximize $\phi^{\rm (w)}$ in \eqref{eq:quadratic}, (OSDEA$(p)$)$^{\rm (w)}$, can be formulated in the same fashion as (OSDEA$(p)$).

The Output Selection for DEA problem where the goal is to maximize $\phi^{\rm (q)}$ in \eqref{eq:quadratic2}, (OSDEA$(p)$)$^{\rm (q)}$, can be formulated similarly to (OSDEA$(p)$). While the feasible region remains the same, the objective function becomes quadratic and the goal is to minimize it, yielding a Mixed Integer Quadratic Programming formulation.

The Output Selection for DEA problem where the goal is to maximize the minimum efficiency $\phi^{\rm (m)}$ in \eqref{eq:quadratic3}, (OSDEA$(p)$)$^{\rm (m)}$ can be written as an MILP. Here, we need to define a new variable $\lambda$ to rewrite the minimum in the objective function, and include the corresponding constraints to ensure that the new variable is well defined.
\begin{equation}
\label{eq:rewritemin}
\lambda \le \sum_{o=1}^O \beta_o^{(k)} y_o^{(k)} \quad k=1,\ldots,K,
\end{equation}
and thus
\begin{eqnarray}
 \mbox{maximize}_{(\bm{\alpha},\bm{\beta},\mathbf{z},\lambda)} \ \displaystyle \lambda
\end{eqnarray}
\noindent {\em s.t.} \hfill (OSDEA$(p)$)$^{\rm (m)}$
 \begin{eqnarray*}
  \displaystyle
  \eqref{eq:DEA1single}-\eqref{eq:DEA4single}; \eqref{eq:DEA5OS}-\eqref{eq:DEA7OS}; \eqref{eq:rewritemin}. && \nonumber
 \end{eqnarray*}

The Output Selection for DEA problem where the goal is to maximize the $\pi$-percentile $\phi^{(\pi)}$ in \eqref{eq:percentile}, (OSDEA$(p)$)$^{(\pi)}$, can also be written as an MILP, similarly as in \cite{benati2015using}. We need to define a new variable $\lambda$ that is equal to the percentile, as well as include the corresponding constraints to ensure that the new variable is well defined. We also need a new binary variable, $\delta^{(k)}$, that is equal to $1$ if the efficiency of DMU $k$,
$\sum_{o=1}^O\beta_o^{(k)} y_o^{(k)}$, is at least $\lambda$ and $0$ otherwise.
\begin{eqnarray}
 \mbox{maximize}_{(\bm{\alpha},\bm{\beta},\mathbf{z},\lambda,\delta)} \ \displaystyle \lambda
\end{eqnarray}
\noindent {\em s.t.} \hfill (OSDEA$(p)$)$^{(\pi)}$
 \begin{eqnarray}
  \displaystyle
  \eqref{eq:DEA1single}-\eqref{eq:DEA4single}; \eqref{eq:DEA5OS}-\eqref{eq:DEA7OS} && \nonumber\\
    \sum_{o=1}^O \beta_o^{(k)} y_o^{(k)} \geq \lambda - M' (1- \delta^{(k)}) && \forall k=1,\ldots,K \label{eq:DEA8singleOSpercentile}\\
    \sum_{k=1}^{K} \delta^{(k)}  = \lfloor K \, \dfrac{\pi}{100}\rfloor \label{eq:DEA9singleOSpercentile}\\
    \delta^{(k)} \in \{0,1\} && \forall  k=1,\ldots,K \label{eq:DEA10singleOSpercentile},
 \end{eqnarray}
with $M'$ a big constant.

\section{Numerical section}\label{sec:numerical}

\begin{figure}[t]
\centering
\includegraphics[scale=0.60]{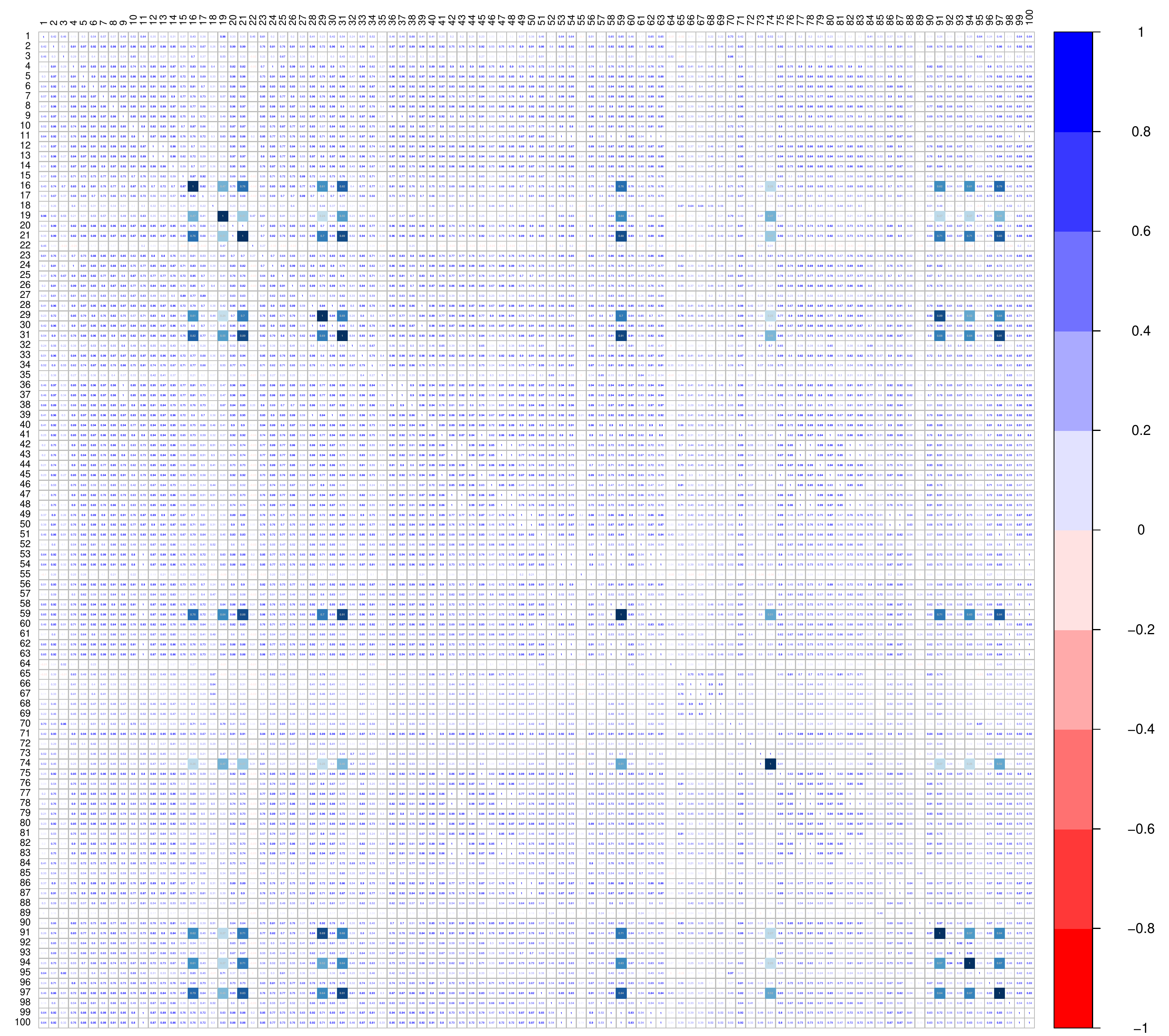}
\caption{Correlation matrix for the outputs, highlighting the correlation between with the selected outputs for $p=10$
\label{fig:Corrp10}}
\end{figure}

In this section, we illustrate the models in previous sections using a real-world dataset in benchmarking of electricity DSOs \cite{agrell2017regulatory,agrell2018theory}. Here, we have $K=182$ DMUs, $O=100$ outputs, and $I=1$ input. As customary, each output has been normalized dividing it by the difference between the maximum and the minimum values of the output.
Figure \ref{fig:Corrp10} displays the correlations between the outputs, with darker colours pointing at higher correlations. This matrix reveals subsets of outputs highly correlated with each other, such as outputs 23 to 31, where correlations are above 0.5, except for $\mbox{corr}(23,27)=0.37$ and $\mbox{corr}(27,29)=0.35$.

The experiments were run on a computer with an Intel$^\circledR$ Core{$^\textrm{TM}$} i7-6700 processor at $3.4$ GHz using $16$ GB of RAM, running Windows 10 Home. All the optimization problems have been solved using Python 3.5 interface \cite{pthn} with Gurobi 7.0.1 solver, \cite{gurobi}.

\begin{figure}[t]
\centering
\includegraphics[scale=0.7]{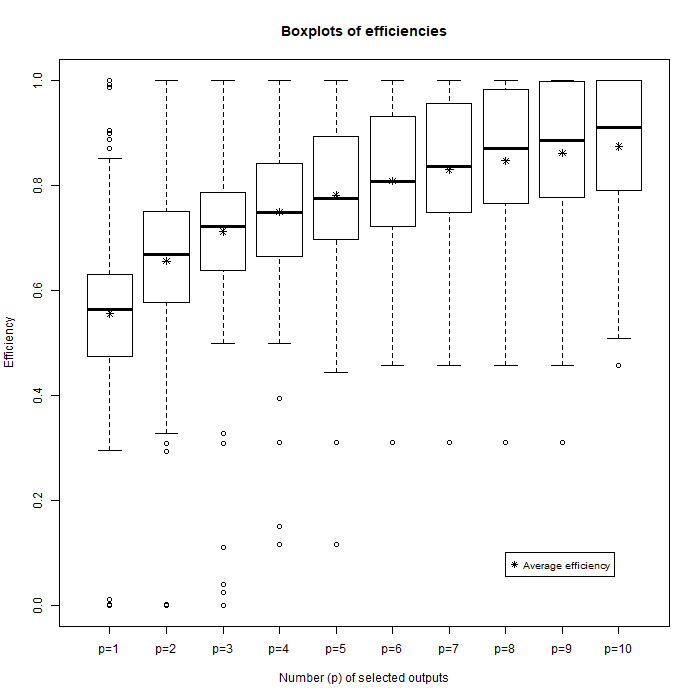}
\caption{Box-and-whiskers plots of efficiencies, including average efficiency, for $p=1,\ldots,10$}
\label{fig:distFSBoxplot}
\end{figure}

\begin{figure}[t]
\centering
\includegraphics[scale=0.70]{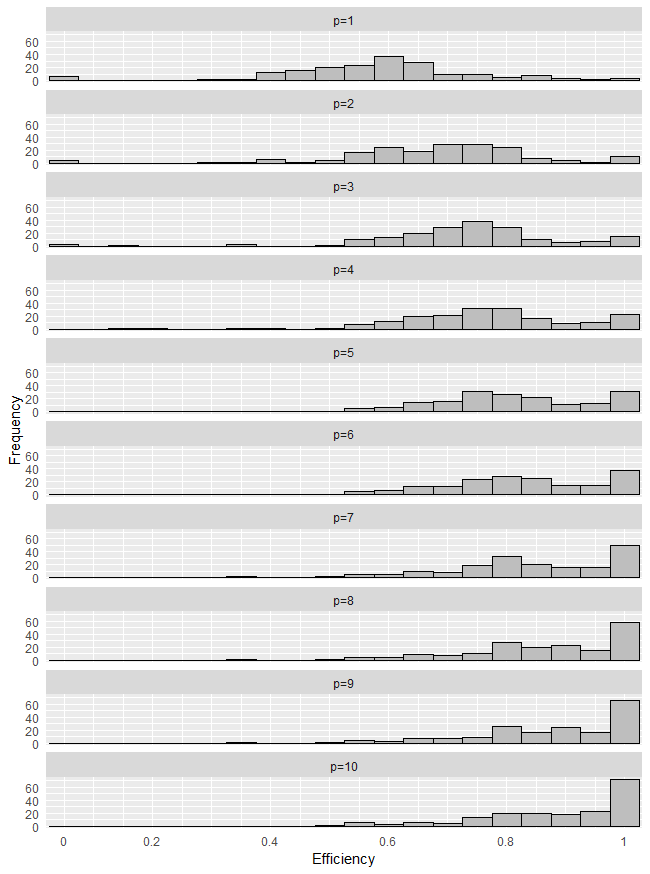}
\caption{Histograms of the distribution of the efficiencies, $p=1,\ldots,10$
\label{fig:distFSallp}}
\end{figure}

\begin{figure}[t]
\centering
\includegraphics[scale=0.90]{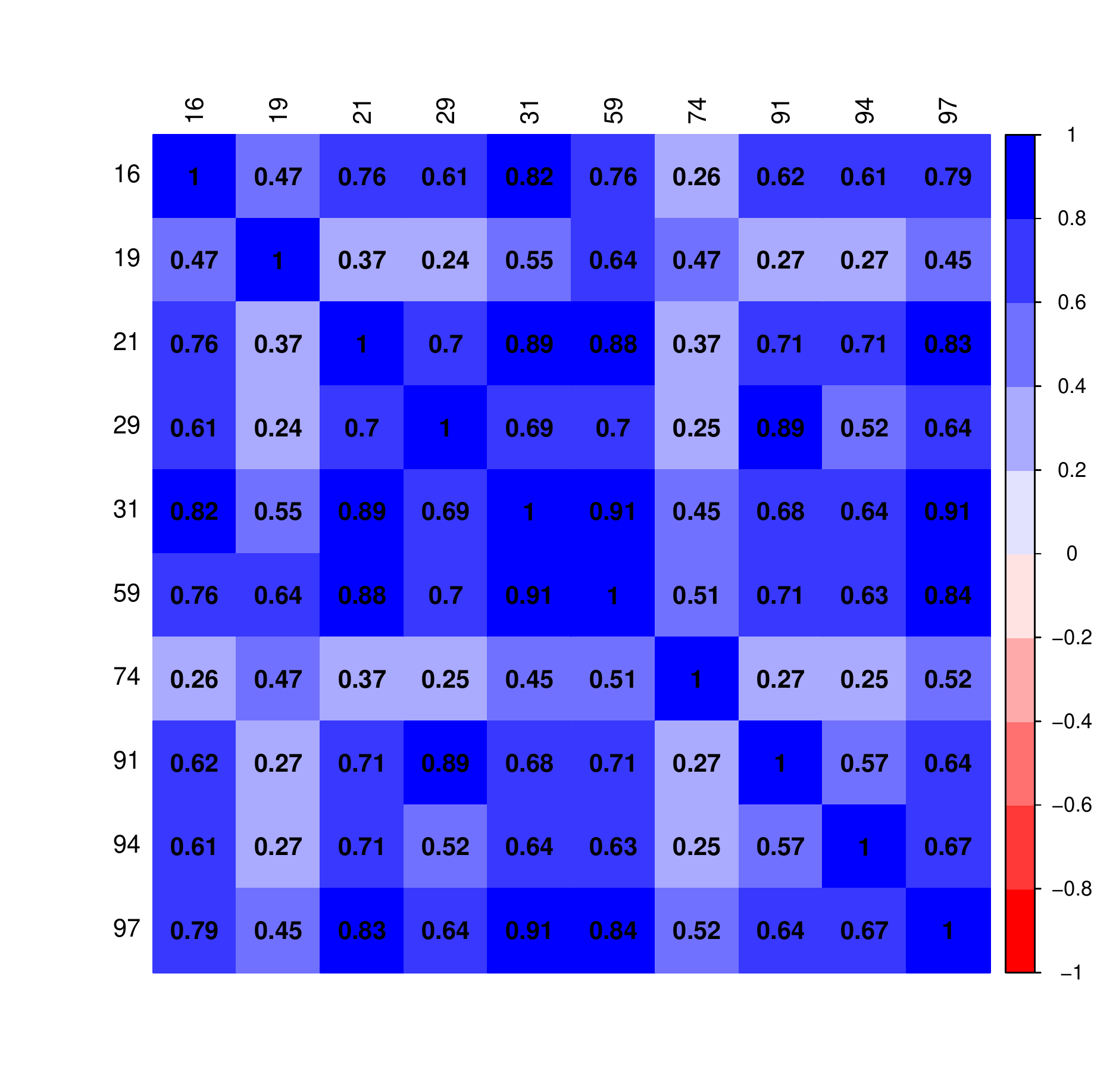}
\caption{Correlation matrix for the selected outputs for $p=10$
\label{fig:Corrp10-selected}}
\end{figure}

We have solved (OSDEA$(p)$) for $p=1,\ldots,10$, with $M$ equal to $1000$. We have run the approach in \cite{pastor2002statistical} to provide (OSDEA$(p)$) with an initial solution. A time limit of 300 seconds has been imposed, although this is not binding for small values of $p$. Once (OSDEA$(p)$) has selected the $p$ outputs, $p=1,\ldots,10$, we calculate the efficiencies of the DSOs obtained with the chosen outputs. The results are summarized in Table \ref{table:summary statistics}, and Figures \ref{fig:distFSallp} and \ref{fig:distFSBoxplot}, while the correlation matrix in Figure \ref{fig:Corrp10-selected} highlights the correlations between the selected outputs.

Table \ref{table:summary statistics} presents summary statistics of the distribution of the efficiencies, namely the minimum, the maximum, the average (i.e., $v(p)$), the standard deviation, the quartiles $q_i$, and the interquartile range (i.e., $q_3$-$q_1$). The last column of this table reports the selected outputs. Figure \ref{fig:distFSBoxplot} displays the box-and-whiskers plots as well as the average efficiency $v(p)$, and Figure \ref{fig:distFSallp} the histograms of the distribution of the efficiencies. The average efficiency improves with the number of selected outputs, $p$, increasing from $0.5555$ to $0.8732$. Figure \ref{fig:distFSBoxplot} shows that the marginal effect of increasing $p$ to $p+1$ is decreasing for this dataset. When looking at the quartiles, we can see that there is a substantial improvement too by increasing $p$. When the number of selected features, $p$, is small the chosen features give poor efficiencies to some of the DMUs. Indeed, for $p\le5$, the minimum efficiency is below $0.1200$. As $p$ increases, we can see that this minimum increases rapidly with $p$. The first quartile increases from $0.4743$ to $0.7902$. Similarly, for the median, we have $0.5637$ to $0.9090$, while for the third quartile $0.6300$ to $1.0000$. The maximum efficiency is already maximal, i.e., equal to $1.0000$, for the smallest value of $p$, and remains like that for all values of $p$ tested. In general, the standard deviation of the efficiencies decreases with $p$, while the interquartile range increases. In terms of the correlations, for small values of $p$, the selected outputs are highly correlated with each other. For instance, for $p=2$, we choose outputs 11 and 31, for which $\mbox{corr}(11,31)=0.91$; while for $p=3$, we choose outputs 21, 31, 54, with $\mbox{corr}(21,31)=0.89$, $\mbox{corr}(21,54)=0.88$ and $\mbox{corr}(31,54)=0.91.$ Actually, for $p=2,\ldots,5$, the smallest correlation between the selected outputs is equal to $0.61$. For $p=10$, there are only two outputs, for which we can find correlations below 0.5, namely outputs 19 and 74.

\begin{table}
\begin{scriptsize}
$$
\begin{tabular}{|c|cc|cc|cccc|l|}
 \hline
$p$&$\min$&$\max$&\textrm{mean}&\textrm{st. dev.}&\textrm{$q_1$}&\textrm{$q_2$}&\textrm{$q_3$}&\textrm{$q_3$-$q_1$}&\textrm{selected features}\\
\hline
1&0.0000&1.0000&0.5555&0.1695&0.4743&0.5637&0.6300&0.1557&\textrm{59}\\
2&0.0006&1.0000&0.6553&0.1708&0.5772&0.6682&0.7482&0.1710&\textrm{11 31}\\
3&0.0009&1.0000&0.7118&0.1643&0.6391&0.7222&0.7839&0.1448&\textrm{21 31 54}\\
4&0.1161&1.0000&0.7487&0.1511&0.6645&0.7479&0.8404&0.1759&\textrm{16 21 31 59}\\
5&0.1161&1.0000&0.7812&0.1494&0.6962&0.7738&0.8911&0.1949&\textrm{16 21 31 59 94}\\
6&0.3105&1.0000&0.8082&0.1402&0.7222&0.8068&0.9305&0.2083&\textrm{16 19 21 31 59 94}\\
7&0.3105&1.0000&0.8290&0.1404&0.7474&0.8355&0.9545&0.2071&\textrm{16 19 21 31 59 91 94}\\
8&0.3105&1.0000&0.8462&0.1402&0.7658&0.8689&0.9802&0.2144&\textrm{16 19 21 31 59 91 94 97}\\
9&0.3105&1.0000&0.8610&0.1370&0.7789&0.8841&0.9972&0.2183&\textrm{16 19 21 31 59 74 91 94 97}\\
10&0.4576&1.0000&0.8732&0.1304&0.7902&0.9090&1.0000&0.2098&\textrm{16 19 21 29 31 59 74 91 94 97}\\
\hline
\end{tabular}
$$
\caption{Summary statistics for the distribution of efficiencies, for $p=1,\ldots,10$}
\end{scriptsize}
\end{table}\label{table:summary statistics}

In real-life applications, the DMUs may be consulted on the chosen outputs. This may be useful to ensure that the resulting choices make conceptual sense. However, such involvement is likely to lead strategic behavior. The evaluated DMUs may try to influence the choice of outputs in their own advantage.
This may lead to a game between the DMUs (the players), each of which has preferred selections of outputs (strategies), and the modeler, who is trying to make a reasonable selection based on the resulting efficiencies of all DMUs (the outcome).

To start investigating the challenges of such strategic behavior, we will consider a simple game or social choice problem. Without loss of generality, we assume that $p$ is fixed. We assume that there are $K$ players, the DMUs, and $K+1$ strategies or choices, namely the selection of outputs $z^{k}(p)$ according to the individual preferences as determined by (OSDEA$^{(k)}(p)$), $k=1,…,K$, and to the joint selection $z(p)$ as determined by (OSDEA$(p)$).
We will think of the joint selection $z(p)$ as the default or status quo selection and the question is now if one of the individual selections $z^{k}(p)$ would be preferred by a large group of DMUs. If so, the modeler is likely to face strong opposition to his proposed selection of outputs.

For DMU $k'$, the selection of outputs made with (OSDEA$^{(k')}(p)$), $z^{(k')}(p)$, is at least as attractive as the one made with (OSDEA$(p)$), $z(p)$, or in other words, the reported efficiency with $z^{(k')}(p)$ is at least as high as the one with $z(p)$. However, for any other DMU $k\not = k'$, it is not clear whether the selection $z^{(k')}(p)$ reports a higher efficiency for DMU $k$ or not. One would like to know how many DMUs prefer the joint strategy $z(p)$ over the $K$ individual strategies $z^{(k')}(p)$. The so-called cross-efficiency measures this preference. Let $\Delta^{(k)}(k')=E^{(k)}(\mathbf{z}^{(k')}(p))-E^{(k)}(\mathbf{z}(p))$, with $k,k'=1,\ldots,K$, be equal to the difference in reported efficiency for DMU $k$ by the joint selection model and the individual selection model for DMU $k'$. We can define
\[\Pi(k') = \frac{100}{K} \, \mbox{cardinality} (\{k \, : \, \Delta^{(k)}(k')>0\}),\]
i.e., the percentage of DMUs that prefer the individual strategy of DMU $k'$ over the joint one.

\begin{figure}[t]
\centering
\includegraphics[scale=0.85]{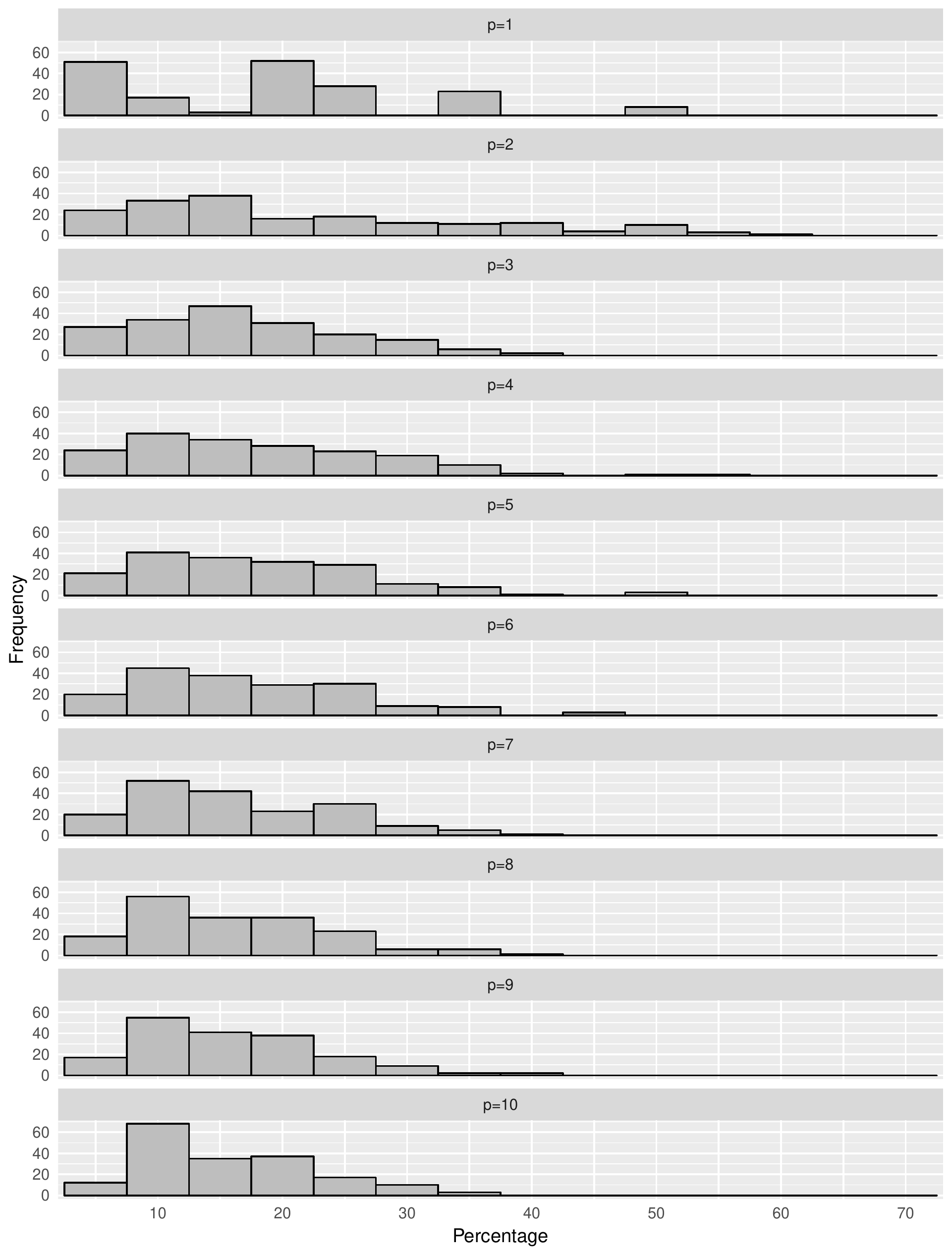}
\caption{Histograms of the distribution of preferences of individual strategies over the joint strategy, $p=1,\ldots,10$
\label{fig:distFSindvsjointallp}}
\end{figure}

In Figure  \ref{fig:distFSindvsjointallp}, we illustrate the share of DMUs that prefer individual selections to the joint selection. As for the joint strategy, a time limit of 300 seconds has been imposed to (OSDEA$^{(k')}(p)$), although this is not binding for any value of $p$. We have binned the support to the individual selections in intervals of width 5\%. The height of the first bar indicates how many individual strategies are preferred by [0\%,5\%) of the DMUs, the height of the second bar corresponds to [5\%,10\%) DMUs supporting it, etc. When $p$ increases, we see that less individual strategies are preferred by many DMUs over the joint strategy. Indeed, for values of $p$ above 5, the joint strategy is supported by at least 50\% of the DMUs over any of the individual strategies, while for values of $p$ above 7, this becomes at least 60\%. We can therefore conclude that, in this simple game, as the model gets larger, it becomes less likely that a large group of DMUs will agree on alternative to the modeler’s selection.

\section{Conclusions} \label{sec:conclusions}

In this paper, we have proposed a single-model approach for feature selection in DEA. When the objective is the average efficiency of the DMUs, the problem can be written as an MILP formulation. We have considered other objectives such as the squared distance to the ideal point, where all the DMUs are efficient, yielding a Mixed Integer Quadratic Programming formulation; and we have shown how to enrich the model to allow for situations where different features come with different costs, e.g., related to data collection or data quality, and where features can be grouped and restrictions can be placed on the use of different groups of features in the specific industrial application. Our numerical section illustrates that we can find good solutions in a reasonable amount of time for the case in which the average efficiency is the goal, which boils down to an MILP.

Our approach deviates from previous literature on feature selection in several ways. It is purely based on mathematical programming as opposed to a mixture of statistical and mathematical programming methods, where the desirable properties above can be modeled in a natural way. It works directly with the original features as opposed to dimensionality reduction techniques, which create artificial features. It focuses on the choice of features from a large set of potential candidate features. Finally, it can handle different objective functions to reflect the underlying objective of the modeler and the application context, e.g., the conflicts between different groups of DMUs in the evaluation.


In this paper, we have also introduced an element of game theory. This is relevant since the evaluated parties in applied projects typically will try to influence the feature selection. It is therefore important to think about the conflict between choosing features from a joint and an individual point of view. We have shown how conflicts can be partially analyzed via the cross-efficiency matrix and we have illustrated the conflict between individual and joint perspectives in the numerical application. In the future, it would however be relevant to further explore these issues. One limitation of our analysis above is that we only consider $K$ specific alternatives to the modeler’s joint selection. In theory, there are, of course, many more alternatives. Indeed, any subset of size $p$ from the set of potential outputs $O$ could potentially muster the support of many DMUs against the modeler’s proposal. The strategic analysis of all possible alternatives is likely to become overwhelming. It may therefore be relevant to introduce some restrictions. One idea is to detect relevant clusters of DMUs and make the selection of features tailored to them. By looking at likely interest groups, the game theoretical analysis may be less complex. Groupings could, for example, refer to small versus large, start-up versus well-established, urban versus rural, and investor-owned versus cooperatively owned DMUs. Also, it would be interesting to add constraints to the feature selection model in order to guarantee the support of a number of DMUs, e.g., a majority of DMUs, hereby modeling more general game theoretical aspects, such as the scope for forming coalitions.

\section*{Acknowledgements}{This research has been financed in part by the EC H2020 MSCA RISE NeEDS Project (Grant agreement ID: 822214); the EU COST Action MI-NET (TD 1409); and research projects MTM2015-65915R, Spain, FQM-329, Junta de Andaluc\'{\i}a, these two with EU ERF funds. This support is gratefully acknowledged.}

\bibliographystyle{plain}
\bibliography{FSandDEA_Revision3_SENT}

\begin{thebibliography}{10}

\bibitem{adler2010improving}
N.~Adler and E.~Yazhemsky.
\newblock {Improving discrimination in data envelopment analysis: PCA--DEA or
  variable reduction}.
\newblock {\em European Journal of Operational Research}, 202(1):273--284,
  2010.

\bibitem{agrell2017regulatory}
P.J. Agrell and P.~Bogetoft.
\newblock Regulatory benchmarking: Models, analyses and applications.
\newblock {\em Data Envelopment Analysis Journal}, 3(1--2):49--91, 2017.

\bibitem{agrell2018theory}
P.J. Agrell and P.~Bogetoft.
\newblock Theory, techniques, and applications of regulatory benchmarking and
  productivity analysis.
\newblock In {\em The Oxford Handbook of Productivity Analysis}. 2018.

\bibitem{allen1997weights}
R.~Allen, A.~Athanassopoulos, R.G. Dyson, and E.~Thanassoulis.
\newblock Weights restrictions and value judgements in data envelopment
  analysis: evolution, development and future directions.
\newblock {\em Annals of Operations Research}, 73:13--34, 1997.

\bibitem{benati2015using}
S.~Benati.
\newblock Using medians in portfolio optimization.
\newblock {\em Journal of the Operational Research Society}, 66(5):720--731,
  2015.

\bibitem{Bertsimas2016ORF}
D.~Bertsimas and A.~King.
\newblock {OR Forum - An Algorithmic Approach to Linear Regression}.
\newblock {\em Operations Research}, 64:2--16, 2016.

\bibitem{bogetoft2013performance}
P.~Bogetoft.
\newblock {\em Performance benchmarking: Measuring and managing performance}.
\newblock Springer Science \& Business Media, 2013.

\bibitem{bogetoft2010benchmarking}
P.~Bogetoft and L.~Otto.
\newblock {\em Benchmarking with Dea, Sfa, and R}, volume 157.
\newblock Springer Science \& Business Media, 2010.

\bibitem{CHARNES1978429}
A.~Charnes, W.W. Cooper, and E.~Rhodes.
\newblock Measuring the efficiency of decision making units.
\newblock {\em European Journal of Operational Research}, 2(6):429--444, 1978.

\bibitem{COOK201945}
W.D. Cook, N.~Ramón, J.L. Ruiz, I.~Sirvent, and J.~Zhu.
\newblock {DEA}-based benchmarking for performance evaluation in
  pay-for-performance incentive plans.
\newblock {\em Omega}, 84:45 -- 54, 2019.

\bibitem{cook2014data}
W.D. Cook, K.~Tone, and J.~Zhu.
\newblock Data envelopment analysis: Prior to choosing a model.
\newblock {\em Omega}, 44:1--4, 2014.

\bibitem{emrouznejad2018survey}
A.~Emrouznejad and G.-L. Yang.
\newblock {A survey and analysis of the first 40 years of scholarly literature
  in DEA: 1978--2016}.
\newblock {\em Socio-Economic Planning Sciences}, 61:4--8, 2018.

\bibitem{fernandez2018stepwise}
F.~Fernandez-Palacin, M.A. Lopez-Sanchez, and M.~Mun{\~o}z-M{\'a}rquez.
\newblock {Stepwise selection of variables in DEA using contribution loads}.
\newblock {\em Pesquisa Operacional}, 38(1):31--52, 2018.

\bibitem{golany1989application}
B.~Golany and Y.~Roll.
\newblock {An application procedure for DEA}.
\newblock {\em Omega}, 17(3):237--250, 1989.

\bibitem{GREEN1996461}
R.H. Green, J.R. Doyle, and W.D. Cook.
\newblock Preference voting and project ranking using dea and cross-evaluation.
\newblock {\em European Journal of Operational Research}, 90(3):461--472, 1996.

\bibitem{gurobi}
{Gurobi Optimization, Inc.}
\newblock Gurobi optimizer reference manual, 2016.

\bibitem{caiLASSODEA}
C.-Y.~Lee J.-Y.~Cai.
\newblock {LASSO variable selection techniques in data envelopment analysis}.
\newblock In {\em The 17th Asia Pacific Industrial Engineering and Management
  Systems Conference (APIEMS 2016)}, Taipei, Taiwan, 2016.

\bibitem{jiang2015dearank}
C.~Jiang and W.~Lin.
\newblock {DEARank: a data-envelopment-analysis-based ranking method}.
\newblock {\em Machine Learning}, 101(1--3):415--435, 2015.

\bibitem{joro2015data}
T.~Joro and P.J. Korhonen.
\newblock Data envelopment analysis.
\newblock In {\em Extension of Data Envelopment Analysis with Preference
  Information}, pages 15--26. Springer, 2015.

\bibitem{landete2017robust}
M.~Landete, J.F. Monge, and J.L. Ruiz.
\newblock {Robust DEA efficiency scores: A probabilistic/combinatorial
  approach}.
\newblock {\em Expert Systems with Applications}, 86:145--154, 2017.

\bibitem{LEE2018}
C.-Y. Lee and J.-Y. Cai.
\newblock {LASSO variable selection in data envelopment analysis with small
  datasets}.
\newblock {\em \emph{Forthcoming in} Omega}, 2018.

\bibitem{li2017variable}
Y.~Li, X.~Shi, M.~Yang, and L.~Liang.
\newblock Variable selection in data envelopment analysis via akaike's
  information criteria.
\newblock {\em Annals of Operations Research}, 253(1):453--476, 2017.

\bibitem{li2017dynamic}
Z.~Li, J.~Crook, and G.~Andreeva.
\newblock {Dynamic prediction of financial distress using Malmquist DEA}.
\newblock {\em Expert Systems with Applications}, 80:94--106, 2017.

\bibitem{luo2012input}
Y.~Luo, G.~Bi, and L.~Liang.
\newblock {Input/output indicator selection for DEA efficiency evaluation: An
  empirical study of Chinese commercial banks}.
\newblock {\em Expert Systems with Applications}, 39(1):1118--1123, 2012.

\bibitem{nataraja2011guidelines}
N.R. Nataraja and A.L. Johnson.
\newblock {Guidelines for using variable selection techniques in Data
  Envelopment Analysis}.
\newblock {\em European Journal of Operational Research}, 215(3):662--669,
  2011.

\bibitem{nunamaker1985using}
T.R. Nunamaker.
\newblock Using data envelopment analysis to measure the efficiency of
  non-profit organizations: A critical evaluation.
\newblock {\em Managerial and Decision Economics}, 6(1):50--58, 1985.

\bibitem{pastor2002statistical}
J.T. Pastor, J.L. Ruiz, and I.~Sirvent.
\newblock {A statistical test for nested radial DEA models}.
\newblock {\em Operations Research}, 50(4):728--735, 2002.

\bibitem{petersen2018directional}
N.C. Petersen.
\newblock {Directional Distance Functions in DEA with Optimal Endogenous
  Directions}.
\newblock {\em Operations Research}, 66(4):1068--1085, 2018.

\bibitem{PODINOVSKI2016916}
V.V. Podinovski.
\newblock {Optimal weights in DEA models with weight restrictions}.
\newblock {\em European Journal of Operational Research}, 254(3):916--924,
  2016.

\bibitem{pthn}
{Python Core Team}.
\newblock {Python: A dynamic, open source programming language}.
\newblock Python Software Foundation, 2015.

\bibitem{qin2014joint}
Z.~{Qin} and I.~{Song}.
\newblock {Joint Variable Selection for Data Envelopment Analysis via Group
  Sparsity}.
\newblock {\em ArXiv e-prints arXiv:1402.3740}, 2014.

\bibitem{ramon2010choice}
N.~Ram{\'o}n, J.L. Ruiz, and I.~Sirvent.
\newblock On the choice of weights profiles in cross-efficiency evaluations.
\newblock {\em European Journal of Operational Research}, 207(3):1564--1572,
  2010.

\bibitem{RUIZ20161}
J.L. Ruiz and I.~Sirvent.
\newblock Common benchmarking and ranking of units with {DEA}.
\newblock {\em Omega}, 65:1 -- 9, 2016.

\bibitem{RUIZ2018}
J.L. Ruiz and I.~Sirvent.
\newblock Performance evaluation through {DEA} benchmarking adjusted to goals.
\newblock {\em \emph{Forthcoming in} Omega}, 2018.

\bibitem{doi101002ev1441}
T.R. Sexton, R.H. Silkman, and A.J. Hogan.
\newblock Data envelopment analysis: Critique and extensions.
\newblock {\em New Directions for Program Evaluation}, 1986(32):73--105, 1986.

\bibitem{sirvent2005monte}
I.~Sirvent, J.L. Ruiz, F.~Borr{\'a}s, and J.T. Pastor.
\newblock {A Monte Carlo evaluation of several tests for the selection of
  variables in DEA models}.
\newblock {\em International Journal of Information Technology \& Decision
  Making}, 4(03):325--343, 2005.

\bibitem{SOLEIMANIDAMANEH20095146}
M.~Soleimani-Damaneh and M.~Zarepisheh.
\newblock {Shannon’s entropy for combining the efficiency results of
  different DEA models: Method and application}.
\newblock {\em Expert Systems with Applications}, 36(3, Part 1):5146--5150,
  2009.

\bibitem{wagner2007stepwise}
J.M. Wagner and D.G. Shimshak.
\newblock {Stepwise selection of variables in Data Envelopment Analysis:
  Procedures and managerial perspectives}.
\newblock {\em European Journal of Operational Research}, 180(1):57--67, 2007.

\bibitem{ZHU2018291}
Q.~Zhu, J.~Wu, and M.~Song.
\newblock Efficiency evaluation based on data envelopment analysis in the big
  data context.
\newblock {\em Computers \& Operations Research}, 98:291--300, 2018.

\end{thebibliography}

\end{document}